\documentclass{amsart}

\usepackage[utf8]{inputenc}
\usepackage{amsmath,amssymb,amsthm,units, stmaryrd,stackrel,relsize,bm}
\usepackage{tikz}
\usepackage{mathdots}

\newtheorem{theorem}{Theorem}

\newtheorem{lemma}[theorem]{Lemma}

\newtheorem{corollary}[theorem]{Corollary}
\newtheorem{proposition}[theorem]{Proposition}

\theoremstyle{definition}
\newtheorem{definition}[theorem]{Definition}

\theoremstyle{remark}
\newtheorem{remark}[theorem]{Remark}
\newtheorem{example}[theorem]{Example}

\usepackage[
    colorlinks=true, 
    citecolor=red,
    linkcolor=blue,
    urlcolor=blue]{hyperref}
\newcommand\kp{{\mathsf{KP}}}
\newcommand\kpi{{\mathsf{KPi}}}

\newcommand\Ord{{\mathsf{Ord}}}

\newcommand{\PI}{\boldsymbol\Pi}

\makeatletter
\def\Ddots{\mathinner{\mkern1mu\raise\p@
\vbox{\kern7\p@\hbox{.}}\mkern2mu
\raise4\p@\hbox{.}\mkern2mu\raise7\p@\hbox{.}\mkern1mu}}
\makeatother
\usetikzlibrary{arrows,positioning} 
\tikzset{
    >=stealth',
    punkt/.style={
           rectangle,
           rounded corners,
           draw=black, very thick,
           text width=6.5em,
           minimum height=2em,
           text centered},
    pil/.style={
           ->,
           thick,
           shorten <=2pt,
           shorten >=2pt,}
}

\begin{document}
\title{The Order of Reflection}
\subjclass[2010]{03D60, 03E45 (Primary); 03D40, 03E10.}
\author{J. P. Aguilera}
\address{Institute of Discrete Mathematics and Geometry, Vienna University of Technology. Wiedner Hauptstra{\ss}e 8--10, 1040 Vienna, Austria.}
\email{aguilera@logic.at}
\date{\today}
\maketitle
\begin{abstract}
Extending Aanderaa's classical result that $\pi^1_1 < \sigma^1_1$, we determine the order between any two patterns of iterated $\Sigma^1_1$- and $\Pi^1_1$-reflection on ordinals. We show that this \emph{linear reflection order} is a prewellordering of length $\omega^\omega$.
This requires considering the relationship between linear and some \emph{non-linear} reflection patterns, such as $\sigma^1_1\wedge\pi^1_1$, the pattern of simultaneous $\Sigma^1_1$- and $\Pi^1_1$-reflection.
%Among our results is the fact that the linear reflection order is cofinal in the non-linear reflection order, though the latter is much longer.

The proofs involve linking the lengths of $\alpha$-recursive wellorderings to various forms of stability and reflection properties satisfied by ordinals $\alpha$ within standard and non-standard models of set theory.
\end{abstract}
%\begin{document}
\tableofcontents
\medskip

\section{Introduction}
Let $L_\alpha$ denote the $\alpha$th level of G\"odel's constructible hierarchy, given by $L_0 = \varnothing$, $L_{\alpha+1} =$ all sets definable over $L_\alpha$ with parameters, and $L_\eta = \bigcup_{\alpha<\eta} L_\alpha$ at limit stages. In $\alpha$-recursion theory, one lifts the usual notion of ``computation'' over the natural numbers (or, equivalently, over $L_\omega$) to $L_\alpha$, for sufficiently closed $\alpha$. As became evident from early work by Kreisel, Kripke, Platek, Sacks, Takeuti, and others (see e.g., Simpson \cite{Si78}), facts about recursion on $L_\alpha$ can be translated into facts about recursion on $L_\omega$ in various ways. In particular, the termination of simple inductive definitions of sets of natural numbers is deeply connected with the reflecting structure of $L$ (see e.g., Cenzer \cite{Ce74} or Aczel and Richter \cite{AcRi74}). The purpose of this article is to study the order in which various reflecting properties given in terms of iterated $\Sigma^1_1$- and $\Pi^1_1$-reflection first occur in the constructible hierarchy.

A formula in the language of set theory is $\Sigma^1_1$ if it contains only existential second-order quantifiers (i.e., ranging over classes) followed by arbitrary first-order quantifiers.
An ordinal $\alpha$ is said to be \emph{$\Sigma^1_1$-reflecting} \index{$\Sigma^1_1$-reflecting} if whenever $\phi$ is a $\Sigma^1_1$ formula in the language of set theory and $a_1,\hdots, a_n$ are finitely many elements of $L_\alpha$, then
\[L_\alpha\models\phi(a_1,\hdots, a_n) \text{ implies } \exists \beta<\alpha\,\Big(a_1,\hdots, a_n \in L_\beta\wedge L_\beta\models\phi(a_1,\hdots, a_n)\Big).\]
Given a class of ordinals $X$, an ordinal $\alpha$ is said to be $\Sigma^1_1$-reflecting \emph{on $X$} if one can additionally demand that the ordinal $\beta$ above belong to $X$. 
The least $\Sigma^1_1$-reflecting ordinal is denoted by $\sigma^1_1$, and $\pi^1_1$ is defined dually. 

An ordinal $\alpha$ is said to be \emph{$\beta$-stable} if $L_\alpha$ is a $\Sigma_1$-elementary substructure of $L_\beta$; in symbols: \index{stability}
$$L_\alpha \prec_1 L_{\beta}.$$
Given an ordinal $\alpha$, write $\alpha^+$ for the smallest admissible ordinal greater than $\alpha$. 
Aczel and Richter \cite{AcRi74} showed that $\pi^1_1 \neq \sigma^1_1$ and that a countable ordinal $\alpha$ is $\Pi^1_1$-reflecting if, and only if, it is $\alpha^+$-stable.
Afterwards, Aanderaa \cite{Aa74} showed that $\pi^1_1 < \sigma^1_1$. Gostanian \cite{Go79} showed that $\sigma^1_1$ is smaller than the least $\alpha$ which is $(\alpha^++1)$-reflecting; in fact, he showed that any $\alpha$ which is both $(\alpha^++1)$-stable and locally countable is also $\Sigma^1_1$-reflecting. Later Gostanian and Hrbacek \cite{GoHr79} employed Gostanian's method to give a new proof of Aanderaa's theorem. A third, apparently folklore proof appears in Simpson \cite{Si78}. Aanderaa's theorem is also an immediate consequence of Proposition \ref{PropSigma11sentence} below, although the proof of Proposition \ref{PropSigma11sentence} has a similar flavor to the argument in Simpson \cite{Si78}.

Let us now generalize the definitions of $\sigma^1_1$ and $\pi^1_1$ as follows:
\begin{definition}
The notion of a \emph{reflection pattern} is given inductively: the empty set is a reflection pattern; if $s$ and $t$ are reflection patterns, then so too are $s\wedge t$, $\sigma^1_1(s)$, and $\pi^1_1(s)$.
\end{definition}
We write $\sigma^1_1$ for $\sigma^1_1(\varnothing)$ and $\pi^1_1$ for $\pi^1_1(\varnothing)$. \index{$\sigma^1_1$}

\begin{definition}
A reflection pattern is \emph{linear} if it contains no conjunctions, and \emph{non-linear} otherwise.
\end{definition}

\begin{definition}
An ordinal is said to be $\sigma^1_1(\varnothing)$-reflecting if it is $\Sigma^1_1$-reflecting; it is said to be $\pi^1_1(\varnothing)$-reflecting if it is $\Pi^1_1$-reflecting.
Let $s$ and $t$ be reflection patterns. Inductively, an ordinal $\alpha$ is said to be $\sigma^1_1(s)$-reflecting if it reflects $\Sigma^1_1$ statements onto $s$-reflecting ordinals; it is said to be $\pi^1_1(s)$-reflecting if it reflects $\Pi^1_1$ statements onto $s$-reflecting ordinals; it is said to be $s\wedge t$-reflecting if it is both $s$-reflecting and $t$-reflecting.
\end{definition}
We may alternate between uppercase $\Sigma$ and $\Pi$ and lowercase $\sigma$ and $\pi$ in speaking about patterns of reflection.

The \emph{ordering problem} is: given two reflection patterns $s$ and $t$, determine whether the least $s$-reflecting ordinal is smaller than the least $t$-reflecting ordinal.
We will identify a pattern $s$ with the least $s$-reflecting ordinal. Thus, instances of the ordering problem are e.g., determining whether
$$\sigma^1_1(\sigma^1_1) <\sigma^1_1(\pi^1_1(\sigma^1_1\wedge\pi^1_1))$$
or whether
$$\pi^1_1(\sigma^1_1) < \pi^1_1(\sigma^1_1) \wedge \pi^1_1(\pi^1_1).$$
Other related problems emerge. For instance, one may ask whether $\sigma^1_1(\sigma^1_1)$ is the least $\Sigma^1_1$-reflecting ordinal which is also a  limit of $\Sigma^1_1$-reflecting ordinals.
(Incidentally, the answer to all three questions is ``no.'') 

In this article, we solve the ordering problem for linear patterns of reflection: we exhibit a way of assigning ordinals to linear patterns in a way that respects their ordering; in particular, we show:

\begin{theorem}
The linear order of reflection is a prewellordering of length $\omega^\omega$.
\end{theorem}

The proof requires analyzing the structure of the non-linear, or full, reflection order, to a certain extent. We shall see that all reflection patterns are witnessed for the first time by ordinals between the least $\alpha$ which is $\alpha^+$-stable and the least $\alpha$ which is $(\alpha^++1)$-stable. In addition, we show:

\begin{theorem}
The linear reflection order is cofinal in the full reflection order.
\end{theorem}

This raises the question of whether the full reflection order also has length $\omega^\omega$. This turns out to be false:

\begin{theorem}
The patterns $\pi^1_1$, $\sigma^1_1$, $\sigma^1_1\wedge\pi^1_1$, and $\sigma^1_1(\sigma^1_1)$ have ranks $1$, $\omega$, $\omega^2$, and $\omega^\omega$ in the reflection order, respectively.
\end{theorem}

In the course of proving these theorems, we find various easier results which we believe to be of independent interest; these are labelled ``propositions.''\\

\iffalse
\paragraph{\bf{Theorem}}\ {\textit{The length of the linear reflection order is $\omega^\omega$. In the full reflection order, the reflection patterns $\pi^1_1$, $\sigma^1_1$, $\sigma^1_1\wedge\pi^1_1$, and $\sigma^1_1(\sigma^1_1)$ have ranks $1$, $\omega$, $\omega^2$, and $\omega^\omega$, respectively. The linear patterns are cofinal in the reflection order. All reflection patterns are first witnessed by ordinals between the least $\alpha$ which is $\alpha^+$-stable, and the least $\alpha$ which is $(\alpha^++1)$-stable.}}\\

\paragraph{\bf{Outline}} The article is divided into seven sections. In Section \ref{SectGandy}, we study connections between the degree of stability of an ordinal $\sigma$ and the lengths of $\sigma$-recursive wellorderings. In Section \ref{SectRTT}, we prove various \emph{reflection transfer} results, of the form ``if $\sigma$ is $s$-reflecting, then it is also $t$-reflecting.'' As an application thereof, we show, in Section \ref{SectPatBSS}, that the rank of the pattern $\sigma^1_1(\sigma^1_1)$ in the reflection order is $\omega^\omega$, and, in Section \ref{SectConjunctionFree}, that the length of the conjunction-free part of the reflection order is also $\omega^\omega$. We finish in Section \ref{SectConcludingPatterns} with some remarks on the length of the reflection order and its generalizations.\\
\fi

\paragraph{\bf{Convention}} Even if not mentioned explicitly, every ordinal in this article is assumed to be both countable and locally countable (i.e., for all $\beta <\alpha$, there is a surjection from $\omega$ to $\beta$ in $L_\alpha$). These are the hypotheses for the theorems of Gostanian and Aczel-Richter mentioned above, respectively.  \index{locally countable}

%\paragraph{\bf{Outline}}
%The article is organized as follows:

\section{Stability and Gandy Ordinals}\label{SectGandy}
For an admissible ordinal $\alpha$, write \index{$\delta_\alpha$}
$$\delta_\alpha = \sup\big\{\delta: \delta \text{ is the length of an $\alpha$-recursive wellordering of a subset of $\alpha$} \big\},$$
where a subset of $\alpha$ is said to be \emph{$\alpha$-recursive} \index{$\alpha$-recursive} if it is $\Delta_1$-definable over $L_\alpha$ with parameters. The value of $\delta_\alpha$ remains unchanged if one replaces ``$\alpha$-recursive'' by ``$\alpha$-r.e.'' in the definition. For every admissible $\alpha$, $\delta_\alpha$ is easily seen to be a limit and e.g., additively indecomposable.
We always have $\delta_\alpha\leq\alpha^+$; an ordinal $\alpha$ is \emph{Gandy} if $\delta_\alpha = \alpha^+$. Gostanian \cite{Go79} showed that $\sigma^1_1$ is the smallest ordinal which is not Gandy. In fact, he showed that a locally countable ordinal is not Gandy if, and only if, it is $\Sigma^1_1$-reflecting. Abramson and Sacks \cite{AbSa76} showed that $(\aleph_\omega^L)^+$ is Gandy, so not every Gandy ordinal is locally countable.\\

Since we know what the degree of stability of $\pi^1_1$ is (viz. $(\pi^1_1)^+$),
a possible first question is that of the degree of stability of $\sigma^1_1$. 

\begin{proposition}\label{TheoremStabilityS11}
$\sigma^1_1$ is not $(\delta_{\sigma^1_1}+1)$-stable.
\end{proposition}
\proof
Let $\delta = \delta_{\sigma^1_1}$. 
Since $\delta<(\sigma^1_1)^+$, it is not admissible. As we observed before, $\delta$ is a limit ordinal; thus, the failure of admissibility must be due to an instance of collection. Choose some $\Delta_0$ formula $\psi$ such that for some $\vec a \in L_\delta$, $L_\delta\not\models \psi(\vec a)\text{-collection}.$
To see that $\sigma^1_1$ is not $(\delta+1)$-stable, consider the formula $\phi$ in the language of set theory asserting that there are sets $A$, $B$ such that:
\begin{enumerate}
\item \label{NotStableCond1} $A$ and $B$ are transitive sets satisfying V=L, $A$ is admissible, $A \in B$, and there is $\vec a \in B$ such that $B$ does not satisfy $\psi(\vec a)\text{-collection}$;
\item \label{NotStableCond2} for each $(\Ord\cap A)$-recursive linear ordering $R \in B$, either there is an infinite descending sequence $b$ through $R$ with $b\in A$, or there is an ordinal $\beta \in B$ and an isomorphism $f \in B$ from $R$ to $\beta$;
\item \label{NotStableCond3} for each $\beta\in B$, there is an $(\Ord \cap A)$-recursive linear ordering $R \in B$ and an isomorphism $f \in B$ from $R$ to $\beta$.
\end{enumerate}
Notice that $\phi$ is a $\Sigma_1$ formula, since the only unbounded quantifier is the one on $B$. Moreover, it does not hold in $L_{\sigma^1_1}$, for the sets $A$ and $B$ would need to be of the form $L_\alpha$ and $L_\beta$, with $\alpha<\beta<\sigma^1_1$. Conditions \eqref{NotStableCond2} and  \eqref{NotStableCond3} together imply that $\beta = \delta_\alpha$, but Gostanian's characterization of $\sigma^1_1$ then implies $\beta = \alpha^+$, contradicting condition \eqref{NotStableCond1}. Finally, it does hold in $L_{\delta+1}$, as witnessed by $A = L_{\sigma^1_1}$ and $B = L_\delta$. To see that \eqref{NotStableCond2} holds, recall a theorem of Gostanian \cite[Theorem 3.2]{Go79} by which if $\alpha$ is $\Sigma^1_1$-reflecting, then every $\alpha$-recursive linear ordering  which is not a wellordering has an infinite descending sequence in $L_{\alpha}$. Thus, every $\sigma^1_1$-recursive linear ordering $R$ either has an infinite descending sequence in $L_{\sigma^1_1}$, or else is isomorphic to some ordinal $\beta<\delta$. One can construct an isomorphism witnessing this by transfinite recursion: at stage $\gamma<\beta$, one has defined $f\upharpoonright \gamma$ and sets $f(\gamma)$ equal to the $R$-least element not in the range of $f\upharpoonright\gamma$. Since this process takes $\beta$-many stages and $R \in L_{\sigma^1_1+1}$, such an isomorphism belongs to $L_{\sigma^1_1 + \delta}$. Since $\delta$ is additively indecomposable, it belongs to $L_\delta$. The proof that \eqref{NotStableCond3} holds is similar.
\endproof

The proof of Proposition \ref{TheoremStabilityS11} shows:
\begin{corollary}\label{CorTheoremStabilityS11}
Suppose $\sigma$ is $\Sigma^1_1$-reflecting and $(\delta_\sigma+1)$-stable. Then, it is a limit of $\Sigma^1_1$-reflecting ordinals.
\end{corollary}

One cannot improve the conclusion of Proposition \ref{TheoremStabilityS11}---every $\Sigma^1_1$-reflecting ordinal is stable to the supremum of its recursive wellorderings:
\begin{proposition}\label{PropAllSigmaAreStable}
Suppose $\sigma$ is $\Sigma^1_1$-reflecting. Then $\sigma$ is $\delta_\sigma$-stable.
\end{proposition}
\proof
Since $\delta$ is a limit ordinal, it suffices to consider arbitrary $\gamma<\delta$ and show that
$$L_{\sigma}\prec_1 L_\gamma.$$
Let $a \in L_{\sigma}$ and $\phi$ be a $\Sigma_1$-formula such that
$L_\gamma\models\phi(a)$.
Without loss of generality, assume that $a$ is an ordinal.
Let $R$ be a $\sigma$-recursive wellordering of length $\gamma$. In particular, $R$ is $\sigma$-r.e., so there is a $\Sigma_1$ formula $\psi$,  such that for all $x,y \in L_{\sigma}$, 
$$xRy \leftrightarrow L_{\sigma}\models\psi(x,y).$$
Let us assume for notational simplicity that $\psi$ is defined without parameters.
Given an ordinal ${\sigma'}$, let $R_{\sigma'}$ be the binary relation given by
$$xR_{\sigma'} y \leftrightarrow L_{\sigma'}\models\psi(x,y).$$
Since $\psi$ is $\Sigma_1$, we have $R_{\sigma'}\subset R_\tau$ whenever ${\sigma'}\leq\tau\leq\sigma$. In particular, $R_{\sigma'}$ is wellfounded for all $\sigma'<\sigma$.

Because $L_\gamma\models\phi(a)$, there is a subset $A$ of $L_{\sigma}$ such that
\begin{enumerate}
\item $A$ codes a model $(M, E)$ of $\kp + V = L$;
\item $M$ has a largest admissible ordinal $\tau$ and $(\tau, E)$ is isomorphic to $(\sigma, \in)$;
\item there is an ordinal $\beta$ of $M$ and a function $f \in M$ which is an isomorphism between $R_\tau^M$ (i.e., $R_\tau$ computed within $M$) and $\beta$, and $L_\beta^M\models \phi(a)$.
\end{enumerate}
The existence of such an $A$ can be expressed by a set-theoretic $\Sigma^1_1$ formula over $L_{\sigma}$ with parameter $a$ (as well as any other parameters involved in the definition of $R$).

By $\Sigma^1_1$-reflection, there is some ${\sigma'}<\sigma$ and some $A_{\sigma'}\subset L_{\sigma'}$ such that $a \in L_{\sigma'}$ and
\begin{enumerate}
\item \label{YesStableCond1} $A_{\sigma'}$ codes a model $(N, F)$ of $\kp + V = L$;
\item \label{YesStableCond3} $N$ has a largest admissible ordinal $\eta$ and $(\eta, F)$ is isomorphic to $({\sigma'},\in)$;
\item \label{YesStableCond4} there is an ordinal $b$ of $N$ and a function $g \in N$ which is an isomorphism between $R_\eta^N$  and $b$, and $L_b^N\models \phi(a)$.
\end{enumerate}
Here and for the rest of our lives, let us identify the wellfounded part of $N$ with its transitive collapse.
Condition \eqref{YesStableCond3}  implies that $L_{\sigma'} \in N$.
By \eqref{YesStableCond4}, there is an ordinal $b$ of $N$ and an isomorphism $g \in N$ from $R_{\sigma'}$ to $b$. 
Because $R_{\sigma'} \subset R$, it is wellfounded, and so $b$ really is an ordinal. Now, $L_b\models\phi(a)$, and $N$ has no admissible ordinals above ${\sigma'}$, so $b<({\sigma'})^+<\sigma$. Since $\phi$ is $\Sigma_1$, we conclude that $L_{\sigma}\models\phi(a)$, as was to be shown.
\endproof

We have shown that $\sigma^1_1$ is $\delta_{\sigma^1_1}$-stable and not $(\delta_{\sigma^1_1}+1)$-stable. 
The proof of Proposition \ref{PropAllSigmaAreStable} illustrates how one derives consequences of an ordinal being $\Sigma^1_1$-reflecting. We shall carry out many similar arguments in the future, perhaps omitting some of the details that show up repeatedly.
We note the following consequence of Proposition \ref{PropAllSigmaAreStable}:

\begin{proposition}\label{PropSigma11sentence}
There is a $\Sigma^1_1$-sentence $\phi$ such that for every countable, locally countable $\sigma$, $L_{\sigma}\models\phi$ if, and only if, $\sigma$ is $\Sigma^1_1$-reflecting or $\Pi^1_1$-reflecting.
\end{proposition}
\proof
Let $\phi$ be the sentence that asserts the existence of some $A\subset L_{\sigma}$ coding a model $(M,E)$ of $\kp + V = L$ containing $\sigma$ and such that 
$$M \models L_{\sigma}\prec_1 L_{\delta_\sigma}.$$
Clearly, every $\Pi^1_1$-reflecting ordinal satisfies this sentence, as does every $\Sigma^1_1$-reflecting ordinal, by Proposition \ref{PropAllSigmaAreStable}. 

Suppose that $L_\sigma \models \phi$, as witnessed by $(M,E)$. Suppose moreover that $\sigma$ is not $\Sigma^1_1$-reflecting, so that $\delta_\sigma = \sigma^+$ by Gostanian's characterization. 
Since $\sigma \in M$, a well-known theorem of F. Ville \index{Ville's Theorem} (see e.g., Barwise \cite{Ba75} for a proof) implies that $L_{\sigma^+}\subset M$. 
Given an arbitrary $\beta<\sigma^+$, we then have $\beta \in M$ and $\beta< \delta_\sigma^M$, for otherwise $\delta_\sigma^M<\sigma^+$, which is impossible, since any $\sigma$-recursive wellordering of a subset of $\sigma$ of length $\delta_\sigma^M$ would belong to $M$.
By choice of $M$,
$$M \models L_{\sigma}\prec_1 L_{\delta_\sigma},$$
and so $M\models L_{\sigma}\prec_1 L_{\beta}$. However, being $\Sigma_1$-elementary is absolute, so we really do have $L_{\sigma}\prec_1 L_{\beta}$ and, since $\beta$ was arbitrary, we have $L_{\sigma} \prec_1 L_{\sigma^+}$, so $\sigma$ is $\Pi^1_1$-reflecting.
\endproof

An immediate consequence is Aanderaa's classical result:
\begin{corollary}[Aanderaa]\label{CorAand}
$\pi^1_1<\sigma^1_1$.
\end{corollary}

Corollary \ref{CorAand} holds in a strong form:
\begin{corollary} \label{CorAand2}
$\sigma^1_1$ reflects $\Sigma^1_1$ sentences on $\Pi^1_1$-reflecting ordinals.
\end{corollary}
\proof
Let $\phi$ be the sentence from Proposition \ref{PropSigma11sentence}. Then, if $\psi$ is another $\Sigma^1_1$ sentence, so is the conjunction $\phi\wedge\psi$.
\endproof
Corollary \ref{CorAand2} is not new; it also follows from the proof of Corollary \ref{CorAand} written down in Simpson \cite{Si78}. Our method for analyzing the reflection order is to prove results akin to Corollary \ref{CorAand2}. Now that we know the degree of stability of $\sigma^1_1$, it is natural to ask what the least ordinal $\alpha$ which is $(\delta_\alpha+1)$-stable is. We shall eventually see that it is rather small and in fact smaller than the successor of $\sigma^1_1$ in the reflection order. We finish this section with some related results that will not be used in future sections.

\begin{proposition}\label{PropReflectingOnReflecting}
Suppose $\sigma$ is locally countable and $\Sigma^1_1$-reflecting on $\Sigma^1_1$-reflecting ordinals. Then $\sigma$ is $(\delta_{\sigma}+1)$-stable.
\end{proposition}
\proof
This is similar to the proof of Proposition \ref{PropAllSigmaAreStable}. Again, it is easy to see that $\delta$ is a limit. Let $\gamma = \delta_{\sigma}+1$ and $a \in L_{\sigma}$ be such that $L_\gamma\models\phi(a)$, for some $\Sigma_1$ formula $\phi$. Let $\psi$ be the $\Sigma^1_1$ formula expressing that $\sigma$ is locally countable and there is a set $A\subset L_\sigma$ coding a model $(M,E)$ of $\kp + V = L$ with $\sigma \in M$ and such that
$$M\models ``L_{\delta_\sigma+1}\models \phi(a).\text{''}$$
Then $L_\sigma\models\psi$. By hypothesis, there is a $\Sigma^1_1$-reflecting $\tau<\sigma$ such that $L_{\tau}\models\psi$. Thus, $\tau$ is locally countable and there is a model $(M,E)$ of $\kp + V = L$ with $\tau \in M$ and such that
$$M\models ``L_{\delta_\tau+1}\models \phi(a).\text{''}$$
By Ville's theorem, $L_{\tau^+}\subset M$ and, since $\tau$ is $\Sigma^1_1$-reflecting, $\delta_{\tau} < \tau^+$. Hence, $M$ computes $\delta_\tau$ and $\delta_\tau+1$ correctly and so we really have $L_{\delta_\tau+1}\models \phi(a)$. Since $\phi$ is $\Sigma_1$, we conclude $L_{\sigma}\models\phi(a)$, as desired.
\endproof

The preceding proof shows that if $\sigma$ is as in Proposition \ref{PropReflectingOnReflecting}, then $\sigma$ is $(\delta_\sigma+2)$-stable, $(\delta_\sigma)^\omega$-stable, etc. It shows that if $f$ is a function on ordinals which is uniformly $\Sigma_1$-definable (with parameters in $L_\sigma$) on e.g., multiplicatively indecomposable levels of $L$ containing all parameters, then $\sigma$ is $f(\delta_\sigma)$-stable.

\begin{definition}
We denote by $\sigma^{1,\ell}_1$ the least $\Sigma^1_1$-reflecting ordinal which is a limit of $\Sigma^1_1$-reflecting ordinals.
\end{definition}

\begin{proposition}\label{PropSigmaLimitDelta+1}
$\sigma^{1,\ell}_1$ is smaller than the least $\alpha$ which is $(\delta_\alpha+1)$-stable.
\end{proposition}
\proof
Let $\alpha$ be as in the statement. We claim that $\alpha$ is $\Sigma^1_1$-reflecting. Otherwise, $\delta_\alpha = \alpha^+$, and so $\alpha$ is $(\alpha^++1)$-stable. But surely $\alpha$ is locally countable, and thus $\Sigma^1_1$-reflecting, by Gostanian's result mentioned in the introduction.
Thus, $\alpha$ is $\Sigma^1_1$-reflecting. By Corollary \ref{CorTheoremStabilityS11}, $\alpha$ is a limit of $\Sigma^1_1$-reflecting ordinals. However, this is expressible in $L_\alpha$; thus, the proof of Proposition \ref{TheoremStabilityS11} shows that $\alpha$ is a limit of ordinals which are both $\Sigma^1_1$-reflecting and limits of $\Sigma^1_1$-reflecting ordinals.
\endproof

As a consequence, we obtain a negative answer to one of the questions posed in the introduction. 

\begin{corollary}
$\sigma^{1,\ell}_1$ does not reflect $\Sigma^1_1$ statements onto $\Sigma^1_1$-reflecting ordinals.
\end{corollary}

We state without proof a result implying that $\sigma^{1,\ell}_1<\pi^1_1(\sigma^1_1)$. Its proof is similar to that of the more powerful Theorem \ref{secondmainwords} below.

\begin{proposition}\label{PropSigmaLimitDelta+2}
For every $\alpha<\pi^1_1(\sigma^1_1)$, there is some $\sigma<\pi^1_1(\sigma^1_1)$ which is both $\Sigma^1_1$-reflecting and $(\delta_\sigma+\alpha)$-stable. 
\end{proposition}

Figure \ref{figureordinalsBelowPS} summarizes the relationships between the ordinals considered so far. We shall also see that 
$$\pi^1_1(\sigma^1_1) < \sigma^1_1(\sigma^1_1).$$

\begin{figure}
\begin{tikzpicture}[node distance=0.3cm, auto,]

\node (ps) {$ \pi^1_1(\sigma^1_1)$};
\node[left=of ps] (dotsagain) {$\hdots$}
edge[->] (ps);
\node[left=of dotsagain] (deltaplustwo) {$L_\alpha\prec_1 L_{\delta_\alpha+2}$}
edge[->] (dotsagain);
\node[left=of deltaplustwo] (deltaplusone) {$L_\alpha\prec_1 L_{\delta_\alpha+1}$}
edge[->] (deltaplustwo);
\node[left=of deltaplusone] (slimit) {$\sigma^{1,\ell}_1$}
edge[->] (deltaplusone);
\node[left=of slimit] (sone) {$\sigma^{1}_1$}
edge[->] (slimit);
\node[left=of sone] (pluszero) {$L_\alpha\prec_1 L_{\alpha^+}$}
 edge[->] (sone);
\end{tikzpicture}
\caption{Ordinals below $\pi^1_1(\sigma^1_1)$} \label{figureordinalsBelowPS}
\end{figure}

%
%
%
\iffalse
\begin{figure}
\begin{tikzpicture}[node distance=0.5cm, auto,]
 %nodes
 \node (plusone) {$L_\alpha\prec_1 L_{\alpha^+ +1}$};
 \node[above=of plusone] (dots) {$\vdots$}
edge[-] (plusone);
%
\node[above=of dots] (ssandps) {$\sigma^1_1(\sigma^1_1(\sigma^1_1))$}
edge[-] (dots);
%
\node[above=of ssandps] (sss) {$\pi^1_1(\sigma^1_1(\sigma^1_1)\wedge \pi^1_1)$}
edge[-] (ssandps);
%
\node[above=of sss] (ssandp) {$\sigma^1_1(\sigma^1_1)\wedge\pi^1_1$}
edge[-] (sss);

%
\node[above=of ssandp] (pss) {$\pi^1_1(\sigma^1_1(\sigma^1_1))$}
edge[-] (ssandp);
%
\node[above=of pss] (ss) {$\sigma^1_1(\sigma^1_1)$}
edge[-] (pss);
%

\node[above=of ss] (ssbis) {$\pi^1_1(\sigma^1_1\wedge\pi^1_1)$}
edge[-] (ss);
%
\node[above=of ssbis] (sandp) {$ \sigma^1_1\wedge\pi^1_1$}
edge[-] (ssbis);
%
\node[above=of sandp] (ps) {$ \pi^1_1(\sigma^1_1)$}
edge[-] (sandp);
%
\end{tikzpicture}
\caption{Some ordinals in normal form} \label{figureordinals}
\end{figure}
%
%
%
\clearpage
\fi
\section{Reflection Transfer Theorems} \label{SectRTT}
In this section, we will present some results on the transfer of reflection properties, i.e., results of the form
\[\text{if }\sigma \text{ is $s$-reflecting,} \text{ then it is $t$-reflecting},\]
where $s$ and $t$ are reflection patterns. Recall our convention that every ordinal considered is countable and locally countable. 
\iffalse
In order to study the reflection order, we focus on ordinals that are small enough. This is not necessary, but will simplify arguments.
\begin{definition}
We say that $\alpha$ is  \emph{subreflecting} if for no $\beta\leq\alpha$ do we have $L_\beta\models$``for every reflection pattern $n$, $(\sigma^1_1)^n$ exists.''\footnote{Replace arbitrary reflection patterns by those of the form $(\sigma^1_1)^n$.}
\end{definition}
Below, we will prove that the sequence 
\[(\sigma^1_1, \sigma^1_1(\sigma^1_1), \sigma^1_1(\sigma^1_1(\sigma^1_1)), \hdots)\]
is cofinal in the reflection order.
Thus, $\alpha$ is subreflecting if, and only if, it is smaller than some reflection ordinal; but the definition given above is more useful, as it is absolute to initial segments of $L$ without strong closure properties. 
\fi
The first five reflection transfer results we present are rather elementary:
\begin{lemma}\label{Propswp<sseasy}
Let $s$ be a reflection pattern and $\sigma$ be an ordinal.
\begin{enumerate}
\item If $\sigma$ is $\Sigma^1_1(\sigma^1_1(s))$-reflecting, then $\sigma$ is $\Sigma^1_1(s)$-reflecting.
\item If $\sigma$ is $\Pi^1_1(\pi^1_1(s))$-reflecting, then $\sigma$ is $\Pi^1_1(s)$-reflecting.
\end{enumerate}
\end{lemma}
\proof
If $\sigma$ is $\Sigma^1_1(\sigma^1_1(s))$-reflecting and $L_\sigma$ satisfies a $\Sigma^1_1$ sentence $\phi$, then, by definition, there is a $\Sigma^1_1(s)$-reflecting $\tau<\sigma$ such that $L_\tau \models\phi$. By $\Sigma^1_1(s)$-reflection, there is an $s$-reflecting $\eta<\tau$ such that $L_\tau\models\phi$. Hence, $\sigma$ is $\Sigma^1_1(s)$-reflecting. The argument for $\Pi^1_1(\pi^1_1(s))$-reflection is similar.
\endproof

\begin{lemma}\label{Propswp<ss}
Let $s$ and $t$ be reflection patterns and $\sigma$ be an ordinal.
\begin{enumerate}
\item \label{Propswp<ss1} If $\sigma$ is $\sigma^1_1(s)\wedge\pi^1_1(t)$-reflecting, then it is $\pi^1_1(\sigma^1_1(s)\wedge t)$-reflecting.
\item \label{Propswp<ss2} If $\sigma$ is $\sigma^1_1(s)\wedge\pi^1_1(t)$-reflecting, then it is $\sigma^1_1(s\wedge \pi^1_1(t))$-reflecting.
\item \label{Propswp<ss3} If $\sigma$ is $\sigma^1_1(s)$-reflecting, then it is $\sigma^1_1(s\wedge \pi^1_1)$-reflecting.
\end{enumerate}
\end{lemma}
\proof
Recall that if an ordinal $\sigma$ is $s$-reflecting, for any nontrivial reflection pattern $s$, then it is recursively inaccessible and, in fact, a limit of recursively inaccessible ordinals. \eqref{Propswp<ss1} then follows from the simple observation that being $\Sigma^1_1(s)$-reflecting is expressible by a $\Pi^1_1$ sentence $\psi$. Thus, if $\sigma$ is $\sigma^1_1(s)\wedge\pi^1_1(t)$-reflecting and satisfies some $\Pi^1_1$-sentence $\phi$, then the conjunction $\phi\wedge\psi$ is also $\Pi^1_1$, and any ordinal satisfying it must be $\Sigma^1_1(s)$-reflecting.

\eqref{Propswp<ss2} is similar. For \eqref{Propswp<ss3}, there are two cases: if $\sigma$ is $\Pi^1_1$-reflecting, then the result follows from \eqref{Propswp<ss2}. If $\sigma$ is not $\Pi^1_1$-reflecting, it is not $\sigma^+$-stable.
Hence, there is a least $\gamma<\sigma^+$ such that $\sigma$ is not $\gamma$-stable, i.e., there is a $\Sigma_1$-formula $\psi$ and some parameter $\alpha<\sigma$ such that $L_\gamma\models\psi(\alpha)$, but $L_\sigma\not\models\psi(\alpha)$.  The remainder of the proof is an adaptation of the proof of Corollary \ref{CorAand} presented in Simpson \cite{Si78}:

Let $\phi$ be the $\Sigma^1_1$-statement expressing that there is a model $(M,E)$ of $\kp + V = L$ end-extending $L_{\sigma+1}$ such that for some $\gamma' \in M$ with $\gamma'<\sigma^{+M}$, $L_{\gamma'}^M\models\psi(\alpha)$ and, moreover, if $\gamma'$ is least such, then $M\models$``$\sigma$ is $({<}\gamma')$-stable.'' Then $L_\sigma\models\phi$. By choice of $\sigma$, there is an $s$-reflecting ordinal $\tau<\sigma$ such that $L_\tau\models\phi$. 
This means that there is a model $(M,E)$ of $\kp + V = L$ end-extending $L_{\tau+1}$ such that for some $\gamma'\in M$ with $\gamma' < \tau^{+M}$, $L^M_{\gamma'}\models\psi(\alpha)$ and, for the least such $\gamma'$, we have $M\models$``$\tau$ is $({<}\gamma')$-stable.''. Since $\psi$ is $\Sigma_1$ and $L_\sigma\not\models\psi(\alpha)$, $\gamma'$ must belong to the illfounded part of $M$. So $\tau$ is $s$-reflecting and, as in the proof of Proposition \ref{PropSigma11sentence}, $\tau$ is $\Pi^1_1$-reflecting. By taking conjunctions as before, one sees that every $\Sigma^1_1$ satisfied by $L_\sigma$ is satisfied by some $(s\wedge\pi^1_1)$-reflecting $\tau<\sigma$, as was to be shown.
\endproof

\begin{example} \label{ExampleSPS<SS}
We claim that
\[\sigma^1_1(\pi^1_1(\sigma^1_1)) < \sigma^1_1(\sigma^1_1).\]
To see this, notice that Lemma \ref{Propswp<ss}\eqref{Propswp<ss1} and  \ref{Propswp<ss}\eqref{Propswp<ss2} imply that
\[\sigma^1_1\wedge\pi^1_1 = \sigma^1_1\wedge\pi^1_1(\sigma^1_1) = \sigma^1_1(\pi^1_1(\sigma^1_1))\wedge\pi^1_1 = \sigma^1_1\wedge\pi^1_1(\sigma^1_1(\pi^1_1(\sigma^1_1)))\]
and thus that
\[\sigma^1_1(\pi^1_1(\sigma^1_1)) < \sigma^1_1\wedge\pi^1_1.\]
On the other hand, Lemma \ref{Propswp<ss}\eqref{Propswp<ss3} implies that
\[\sigma^1_1(\sigma^1_1) = \sigma^1_1(\sigma^1_1\wedge\pi^1_1),\]
and so
\[\sigma^1_1(\pi^1_1(\sigma^1_1)) < \sigma^1_1\wedge\pi^1_1 < \sigma^1_1(\sigma^1_1),\]
as claimed. \qed
\end{example}
A natural question is whether one can strengthen $\pi^1_1$ in the statement of Lemma \ref{Propswp<ss}\eqref{Propswp<ss3} and, in particular, whether $\sigma^1_1$ is $\sigma^1_1(\pi^1_1(\pi^1_1))$-reflecting. 
By generalizing the proof of Lemma \ref{Propswp<ss}\eqref{Propswp<ss3}, we see that the answer is ``yes.''

\begin{definition} 
Let $s$ be a reflection pattern. An ordinal $\alpha$ is \emph{$\beta$-stable on $s$} if whenever $L_\beta$ satisfies a $\Sigma_1$ sentence $\phi(L_\alpha)$ with additional parameters in $L_\alpha $, there is an $s$-reflecting $\gamma<\alpha$ such that $L_{\gamma^+}\models \phi(L_\gamma)$. \index{stability!on a reflection pattern}
\end{definition}
We caution the reader that an ordinal $\alpha$ being $\beta$-stable on $\varnothing$ is not the same as it being $\beta$-stable, for the first definition allows $L_\alpha$ as a parameter. We do have the following:

\begin{lemma}\label{LemmaCharStableOnS}
Let $s$ be a reflection pattern. The following are equivalent:
\begin{enumerate}
\item $\alpha$ is $\alpha^+$-stable on $s$;
\item $\alpha$ is $\Pi^1_1(s)$-reflecting.
\end{enumerate}
\end{lemma}
We omit the proof of Lemma \ref{LemmaCharStableOnS}, which is a simple adaptation of Aczel and Richter's characterization of $\Pi^1_1$-reflection.

\begin{theorem}\label{firstmainwords}
Let $s$ be a reflection pattern.
Suppose $\sigma$ is $\sigma^1_1(s)$-reflecting.  Then, it is $\sigma^1_1(s\wedge \pi^1_1(s))$-reflecting.
\end{theorem}
\proof
The conclusion of the theorem follows from Lemma \ref{Propswp<ss} if $\sigma$ is $\pi^1_1(s)$-reflecting, so we may assume that it is not. 

Since $\sigma$ is not $\pi^1_1(s)$-reflecting, it is not $\sigma^+$-stable on $s$, so there is a least $\beta<\sigma^+$ and a $\Sigma_1$-formula $\exists x\,\phi(y,x)$ such that $L_\beta\models\exists x\,\phi(L_\sigma,x)$ and whenever $\gamma<\sigma$ and $\gamma$ is $s$-reflecting, then $L_{\gamma^+} \not\models\exists x\,\phi(L_\gamma,x)$. Let $\psi$ be the formula expressing that there is a model $M$ of $\kp + V = L$ such that
\begin{enumerate}
\item $M$ contains $\sigma$.
\item $M\models$ ``$\exists x\, \phi(L_\sigma,x)$ and, letting $\gamma$ be least such that $M\models\phi(L_\sigma,a)$ for some $a \in L_\gamma$, $\sigma$ is ${<}\gamma$-stable on $s$.''
\end{enumerate}
By reflection, there is $\tau<\sigma$ with $L_\tau\models\psi$, as witnessed by some model $N$ which end-extends $L_{\tau^+}$. Since $\sigma$ is $\sigma^1_1(s)$-reflecting, we may assume that $\tau$ is $s$-reflecting.
Let $\gamma$ be $N$-least such that $N\models\exists x\in L_\gamma\, \phi(L_\tau,x)$. Then, we cannot have $\gamma < \tau^+$, for otherwise $\tau$ is an $s$-reflecting ordinal such that $L_{\tau^+}\models\exists x\, \phi(L_\tau,x)$, contradicting the choice of $\phi$. Thus, $\gamma$ belongs to the illfounded part of $N$ and, in $N$, $\tau$ is ${<}\gamma$-stable on $s$. Since $\tau$ is recursively inaccessible (this can be assumed also if $s = \varnothing$), $N$ is correct about $s$-reflection below $\tau$, so an argument as before shows that $\tau$ is $\tau^+$-stable on $s$ and thus $\Pi^1_1(s)$-reflecting.
\endproof

\begin{example}
By repeatedly applying Theorem \ref{firstmainwords}, we obtain
$$\sigma^1_1 = \sigma^1_1(\pi^1_1) = \sigma^1_1(\pi^1_1(\pi^1_1)) = \sigma^1_1(\pi^1_1(\pi^1_1(\pi^1_1))) =  \dots.$$
This implies the sequence 
of inequalities 
$$\pi^1_1 < \pi^1_1(\pi^1_1) < \pi^1_1(\pi^1_1(\pi^1_1)) < \dots < \sigma^1_1,$$
which strengthens Corollary \ref{CorAand}. 
\end{example}

The following strengthening of Proposition \ref{PropAllSigmaAreStable} is proved similarly:
\begin{lemma}\label{LemmaStrongAllSigmaAreStable}
Suppose $\sigma$ is $\Sigma^1_1(s)$-reflecting. Then, it is $\delta_\sigma$-stable on $s$.
\end{lemma}
\proof
Let $\theta(L_\sigma)$ be a $\Sigma_1$ sentence with parameters in $L_\sigma$, say, of the form $\exists x\,\theta_0(x,L_\sigma)$. Let $\eta<\delta_\sigma$ and $b \in L_\eta$ be such that $L_{\delta_\sigma}\models\theta_0(b,L_\sigma)$. Since $\eta<\delta_\sigma$, there is a $\sigma$-recursive wellorder $R$ of length $\eta$. Let $\psi$ be the sentence asserting the existence of a model $M$ of $\kpi$ such that\footnote{$\kpi$ is the extension of $\kp$ by an axiom asserting that every set is contained in an admissible set.}
\begin{enumerate}
\item $M$ end-extends $L_{\sigma+1}$;
\item in $M$, $R$ is isomorphic to an ordinal $\eta_0$ and there is $b_0 \in L_{\eta_0}^M$ such that $L_{\eta_0}^M \models\theta_0(b_0,L_\sigma)$.
\end{enumerate}
Then $L_\sigma\models\psi$. Moreover, $\psi$ is $\Sigma^1_1$ so, by reflection, there is an $s$-reflecting $\tau<\sigma$ such that $L_\tau\models\psi$, as witnessed by some model $N$ which end-extends $L_{\tau^+}$. Now, in $N$, $L_{\tau^+}^N\models\theta_0(b_0, L_\tau)$ for some $b_0 \in L^N_{\eta_0}$, where $\eta_0$ is some $N$-ordinal isomorphic to $R_\tau$. However, $R_\tau\subset R$, since $R$ is $\sigma$-recursive, and $\theta_0$ is $\Sigma_0$, so we really have $L_{\tau^+}\models\theta(L_\tau)$, as desired.
\endproof

The following theorem, although perhaps odd-looking at first, is crucial for our analysis of the reflection order.
\begin{theorem}
\label{secondmainwords}
Let $s$ be a reflection pattern.
Suppose $\sigma$ is $\Pi^1_1(\sigma^1_1(s))$-reflecting but not $\Sigma^1_1$-reflecting. Then $\sigma$ is $\Pi^1_1(s)$-reflecting.
\end{theorem}
\proof
Suppose $\sigma$ is $\Pi^1_1$-reflecting on $\Sigma^1_1(s)$-reflecting ordinals but not $\Sigma^1_1$-reflecting.
Let $\psi$ be the statement expressing that whenever $(M,E)$ is an end-extension of $L_{\sigma+1}$ satisfying $\kpi$, then $M\models$``$\sigma$ is not $\Sigma^1_1(s)$-reflecting.'' This sentence is $\Pi^1_1$ and thus cannot be satisfied by $L_\sigma$, for otherwise it would be reflected to a $\Sigma^1_1(s)$-reflecting ordinal. But clearly $L_\alpha$ cannot satisfy $\psi$ if $\alpha$ is $\Sigma^1_1(s)$-reflecting.

Thus, $L_\sigma\not\models \psi$, so there is a model $M$ of $\kpi$ end-extending $L_{\sigma+1}$ such that
\begin{align*}
\text{$M\models$``$\sigma$ is $\Sigma^1_1(s)$-reflecting.''}
\end{align*}
%Notice that $M$ contains $\omega$, so all reflection patterns in the sense of $M$ are reflection patterns outside of $M$. Moreover, f
For ordinals $\tau<\sigma$, whether $\tau$ is $t$-reflecting is computed correctly by, say, $\tau^{++}$, and thus too by $M$, for any reflection pattern $t$. By Lemma \ref{LemmaStrongAllSigmaAreStable} applied within $M$, 
\begin{align*}
\text{$M\models$``$\sigma$ is $\delta_\sigma$-stable on $s$.''}
\end{align*}

Let $\phi$ be a $\Pi^1_1$ statement, and $a \in L_\sigma$ be a parameter such that $L_\sigma \models\phi(a)$. By Barwise-Gandy-Moschovakis \cite{BGM}, there is a $\Sigma_1$ formula $\phi^*$ such that for all admissible $\alpha$ with $a \in L_\alpha$, $L_\alpha\models\phi(a)$ if, and only if, $L_{\alpha^+}\models\phi^*(a,L_\alpha)$; thus, $L_{\sigma^+} \models\phi^*(a,L_\sigma)$. Let $b \in L_{\sigma^+}$ be a witness for $\phi^*$ and let $\gamma<\sigma^+$ be large enough so that $b \in L_\gamma$. Since $\sigma$ is not $\Sigma^1_1$-reflecting (in the real world), $\delta_\sigma = \sigma^+$, and thus
$$\gamma < \delta_\sigma^M.$$
Since $\phi^*$ is $\Sigma_1$,
\begin{align*}
\text{$M\models$``$L_{\delta_\sigma}\models \phi^*(a,L_\sigma)$,''}
\end{align*}
so by the $\delta_\sigma$-stability of $\sigma$ on $s$ within $M$, there is an $s$-reflecting $\tau<\sigma$ such that 
\begin{align*}
\text{$M\models$``$L_{\tau^+}\models \phi^*(a,L_\tau)$.''}
\end{align*}
Since $\tau^+ < \sigma$, we really do have 
$$L_{\tau^+}\models \phi^*(a,L_\tau),$$
and so $L_\tau\models\phi(a)$. This completes the proof of the theorem.
\endproof

\begin{remark}
The assumption that $\sigma$ is not $\Sigma^1_1$-reflecting cannot be removed from the statement of Theorem \ref{secondmainwords}. To see this, let $s$ be the trivial pattern. By Lemma \ref{Propswp<ss}\eqref{Propswp<ss1}, $\sigma^1_1\wedge\pi^1_1$ is $\Pi^1_1(\sigma^1_1)$-reflecting. However, $\sigma^1_1\wedge\pi^1_1$ is not $\Pi^1_1(\pi^1_1)$-reflecting, for being $\Sigma^1_1$-reflecting is expressible by a $\Pi^1_1$-formula, and thus every ordinal which is $(\sigma^1_1\wedge\pi^1_1(\pi^1_1))$-reflecting is also
$\pi^1_1(\sigma^1_1\wedge\pi^1_1)$-reflecting and, in particular, a limit of $(\sigma^1_1\wedge\pi^1_1)$-reflecting ordinals.\qed
\end{remark}

\begin{remark}
One cannot improve the statement of Theorem \ref{secondmainwords} to conclude that $\sigma$ is $\Pi^1_1(\sigma^1_1(s)\wedge s)$-reflecting, for let $s = \pi^1_1$. Then, $\pi^1_1(\sigma^1_1(\pi^1_1)) = \pi^1_1(\sigma^1_1)$ by Lemma \ref{Propswp<ss}\eqref{Propswp<ss3}. However, as in Example \ref{ExampleSPS<SS},
\[\pi^1_1(\sigma^1_1) < \sigma^1_1(\pi^1_1(\sigma^1_1)) < \sigma^1_1\wedge\pi^1_1.\]
This is in contrast to Theorem \ref{firstmainwords}.\qed
\end{remark}

\begin{example}
By combining Theorems \ref{secondmainwords} and \ref{firstmainwords}, one sees that 
$$\pi^1_1(\sigma^1_1) = \pi^1_1(\sigma^1_1) \wedge \pi^1_1(\pi^1_1) = \pi^1_1(\sigma^1_1) \wedge \pi^1_1(\pi^1_1(\pi^1_1)) = \hdots$$
Since we have seen that these reflection patterns are all smaller than $\sigma^1_1(\sigma^1_1)$ and $\sigma^1_1\wedge\pi^1_1$, it follows that  $\sigma^1_1$ and $\pi^1_1(\sigma^1_1)$ have order-types $\omega$ and $\omega+1$ in the reflection order, respectively.\qed
\end{example}

\begin{example}
Let us present a proof of the inequality
\[\pi^1_1\Big(\sigma^1_1\wedge\pi^1_1(\sigma^1_1\wedge \pi^1_1)\Big) < \sigma^1_1(\sigma^1_1).\]
First, apply Theorem \ref{firstmainwords} to see that
\[\sigma^1_1(\sigma^1_1) = \sigma^1_1(\sigma^1_1\wedge \pi^1_1(\pi^1_1)).\]
Then, apply Lemma \ref{Propswp<ss} to see that 
$$\sigma^1_1\wedge \pi^1_1(\pi^1_1) = \sigma^1_1\wedge \pi^1_1(\sigma^1_1\wedge\pi^1_1),$$ 
so that
\[\sigma^1_1\Big(\sigma^1_1\wedge \pi^1_1(\pi^1_1)\Big) = \sigma^1_1\Big(\sigma^1_1\wedge \pi^1_1(\sigma^1_1\wedge\pi^1_1)\Big).\]
Finally, by Theorem \ref{firstmainwords}, every $\sigma^1_1(\sigma^1_1\wedge \pi^1_1(\sigma^1_1\wedge\pi^1_1))$-reflecting ordinal is also $\sigma^1_1(\pi^1_1(\sigma^1_1\wedge \pi^1_1(\sigma^1_1\wedge\pi^1_1)))$-reflecting,
so that $\sigma^1_1(\sigma^1_1)$ is a limit of $\pi^1_1(\sigma^1_1\wedge \pi^1_1(\sigma^1_1\wedge\pi^1_1))$-reflecting ordinals. \qed
\end{example}

We finish this section with a final reflection transfer theorem.
It is a strengthening of Theorem \ref{secondmainwords} which clarifies the hypothesis on $\sigma$ not being $\Sigma^1_1$-reflecting. We state it separately, however, since the proof is longer and the result is not used afterwards.

\begin{theorem}\label{fifthmainwords}
Suppose $\sigma$ is $\pi^1_1(\sigma^1_1(t\wedge \pi^1_1(s)))$-reflecting. Then, one of the following holds:
\begin{enumerate}
\item $\sigma$ is $\pi^1_1(s)$-reflecting; or
\item $\sigma$ is $\sigma^1_1(t\wedge \pi^1_1(s))$-reflecting.
\end{enumerate}
\end{theorem}
\proof
Suppose $\sigma$ is $\pi^1_1(\sigma^1_1(t\wedge \pi^1_1(s)))$-reflecting but not $\pi^1_1(s)$-reflecting. Let $\theta$ be a $\Sigma^1_1$ sentence with parameters in $L_\sigma$ such that
\[L_\sigma\models\theta;\]
we need to find a $(t\wedge \pi^1_1(s))$-reflecting $\tau<\sigma$ such that
\[L_\tau\models\theta.\]
By Barwise-Gandy-Moschovakis \cite{BGM}, \index{Barwise-Gandy-Moschovakis} there is a $\Pi_1$ formula $\theta^*(a)$ such that for every admissible $\alpha$ containing the parameters of $\theta$, $L_\alpha\models\theta$ if, and only if, $L_{\alpha^+}\models\theta^*(L_\alpha).$
In particular, 
\[L_{\sigma^+}\models\theta^*(L_\sigma).\]

Since $\sigma$ is not $\pi^1_1(s)$-reflecting, there is a least $\beta<\sigma^+$ such that $\sigma$ is not $\beta$-stable on $s$. Because $\theta^*$ is $\Pi_1$,
\[L_{\beta}\models\theta^*(L_\sigma).\]
Let $\phi$ be the sentence asserting the non-existence of a model $M$ of $\kpi + V = L$ end-extending $L_{\sigma+1}$ in which $\sigma$ is $\sigma^1_1(t\wedge \pi^1_1(s))$-reflecting. This is a $\Pi^1_1$ sentence and thus cannot be satisfied by any $\pi^1_1(\sigma^1_1(t\wedge \pi^1_1(s)))$-reflecting ordinal and, in particular, by $\sigma$. Thus, there is a model $M$ of $\kpi + V = L$ end-extending $L_{\sigma+1}$ and such that
\[M\models \text{``$\sigma$ is $\sigma^1_1(t\wedge \pi^1_1(s))$-reflecting.''}\]
Let $\chi$ be the sentence asserting the existence of a model $N$ of $\kp + V = L$ such that
\begin{enumerate}
\item $N$ contains $\sigma$;
\item in $N$, letting $\beta$ be least such that $\sigma$ is not $\beta$-stable on $s$, we have
\[N\models ``L_\beta\models \theta^*(L_\sigma).\text{''}\]
\end{enumerate}
Since $\beta <\sigma^+$ and $M$ must end-extend $L_{\sigma^+}$, $\beta \in M$ and $M$ is correct about $\beta$ being the least ordinal at which $\sigma$ fails to be stable on $s$.
Thus, we have
\[M\models ``L_\sigma\models \chi,\text{''}\]
as witnessed, say, by $L^M_{\sigma^+}$. Within $M$, $\sigma$ is $\sigma^1_1(t\wedge \pi^1_1(s))$-reflecting and thus there is some $\tau<\sigma$ such that, 
\[M\models ``L_\tau\models \chi,\text{''}\]
and so we really do have 
\[L_\tau\models \chi.\]
Moreover, $M$ is correct about reflection below $\sigma$, so we may assume that $\tau$ is $(t\wedge \pi^1_1(s))$-reflecting. By the definition of $\chi$, there is
a model $N$ of $\kp + V = L$ such that
\begin{enumerate}
\item $N$ contains $\tau$;
\item in $N$, letting $\gamma$ be least such that $\tau$ is not $\gamma$-stable on $s$, we have
\[N\models ``L_\beta\models \theta^*(L_\tau).\text{''}\]
\end{enumerate} 
Since $\tau$ is $\pi^1_1(s)$-reflecting, it is $\tau^+$-stable on $s$, and thus $\gamma$ cannot be a true ordinal smaller than $\tau^+$. By Ville's Theorem, $N$ must end-extend $L_{\tau^+}$. Because $\theta^*$ is $\Pi_1$, it follows that
\[L_{\tau^+} \models \theta^*(L_\tau),\]
and thus, that
\[L_\tau\models\theta,\]
as was to be shown.
\endproof

\section{Patterns Below $\sigma^1_1(\sigma^1_1)$}\label{SectPatBSS}
In this section, we describe the reflection order below $\sigma^1_1(\sigma^1_1)$. The remaining sections do not depend on this one, so the reader who so desires should feel free to skip ahead. To ease notation, we shall omit subscripts and superscripts and simply write $\sigma$ for $\sigma^1_1$ and $\pi$ for $\pi^1_1$. We shall also sometimes omit parentheses; thus, e.g., we will write
$$\sigma\wedge\pi\sigma$$
instead of
$$\sigma^1_1\wedge\pi^1_1(\sigma^1_1).$$
We will also express concatenation by direct juxtaposition, so that e.g., if $s = \sigma\wedge\pi$, then
\[ss = \sigma\wedge\pi(\sigma\wedge\pi).\]
We remind the reader one last time of our convention on all ordinals being countable and locally countable.
 The following notation will be useful:
\begin{definition}
Let $s$ and $t$ be reflection patterns. We write $s \equiv t$ if for every  ordinal $\alpha$, $\alpha$ is $s$-reflecting if, and only if, it is $t$-reflecting.
\end{definition}
\begin{definition}
Let $k \in\mathbb{N}$ and $s$ be a reflection pattern.
We write $c^k_0s = s$; inductively,
\[c_{n+1}^ks = \sigma\wedge\pi^k\sigma\pi c^k_ns.\]
We write $c^k_n$ for $c^k_ns$, where $s$ is the empty pattern.
\end{definition}
We remark that, in particular, $c^0_n = (\sigma\pi)^{n}$.

\begin{lemma}\label{LemmaBelowOmegaCubed}
For every $n,k \in\mathbb{N}$ and every reflection pattern $s$, every $(\sigma\wedge \pi^{k+1}s)$-reflecting ordinal is $(\pi c^k_ns)$-reflecting.
\end{lemma}
\proof
We first show by induction  on $n$ that
\[\sigma\wedge \pi^{k+1} \equiv \sigma \wedge \pi c^{k}_n.\]
Suppose that every $(\sigma\wedge \pi^{k+1})$-reflecting ordinal is $(\pi c^k_n)$-reflecting, i.e., that
$$\sigma\wedge\pi^{k+1} \equiv \sigma\wedge \pi c^k_n \equiv \sigma\wedge \pi (c^k_n\wedge\pi^k).$$
After some applications of Lemma \ref{Propswp<ss}, we have
\begin{align*}
\sigma\wedge \pi c^k_n
&\equiv \sigma\pi c^k_n \wedge \pi c^k_n\\
&\equiv \sigma\pi c^k_n \wedge \pi^{k+1}\\
&\equiv \sigma\pi c^k_n \wedge \pi^{k+1}\sigma\pi c^k_n\\
&\equiv \sigma \wedge \pi^{k+1}\sigma\pi c^k_n\\
&\equiv \sigma \wedge \pi(\pi^{k}\sigma\pi c^k_n)\\
&\equiv \sigma \wedge \pi(\sigma\wedge\pi^k\sigma\pi c^k_n)\\
&\equiv \sigma \wedge \pi c^{k}_{n+1},
\end{align*}
as desired. The argument given also shows that every every $(\sigma\wedge \pi^{k+1}s)$-reflecting ordinal is $(\sigma\wedge\pi c^k_ns)$-reflecting, although it does not prove the converse (which, incidentally, is not true).
\endproof

\begin{corollary}\label{LemmaBelowOmegaCubedCorollary}
For every $n,k \in\mathbb{N}$, every $l\leq k$, and every reflection pattern $s$, every $(\sigma\wedge \pi^{k+1}s)$-reflecting ordinal is $(\pi^{k+1}\sigma\pi c^l_{n}s)$-reflecting.
\end{corollary}
\proof
Suppose $\alpha$ is $(\sigma\wedge \pi^{k+1}s)$-reflecting. By Lemma \ref{LemmaBelowOmegaCubed}, 
\[\text{$\alpha$ is $\pi c^k_{n+1}s$-reflecting.}\]
By definition,
\[\pi c^k_{n+1}s \equiv \pi(\sigma\wedge\pi^k\sigma\pi c^k_n s);\]
in particular,
\[\text{$\alpha$ is $\pi^{k+1}\sigma\pi c^k_{n}s$-reflecting.}\]
We may thus apply Lemma \ref{Propswp<sseasy} $n$ times to see that every $c^k_n s$-reflecting ordinal is also $c^l_n$-reflecting, from which the result follows.
\endproof

\iffalse
\begin{lemma}\label{LemmaSSNFSigmaPre}
Suppose $m_0, m_1, \hdots, m_l$ are natural numbers all smaller than $k$. Then, for any reflection pattern $s$,
\[c^k_1s \equiv c^k_1c^{m_0}_1\hdots c^{m_l}_1\sigma\pi s.\]
\end{lemma}
\proof
We prove it for the case $l = 0$; the general case is similar.
Letting $m = m_0$, we have 
\begin{align*}
c^k_1c^{m}_1s 
&\equiv \sigma\wedge \pi^k\sigma\pi(\sigma\wedge \pi^m\sigma\pi s)\\
&\equiv \sigma\wedge \pi^k\sigma\pi(\sigma\wedge \pi^ms)\\
&\equiv \sigma\wedge \pi^k\sigma\pi(\sigma\pi^ms\wedge \pi^ms).
\end{align*}
By Lemma \ref{LemmaBelowOmegaCubed}, every $(\sigma\wedge\pi^k\sigma\pi s)$-reflecting ordinal is $(\sigma\wedge \pi^kc^m_1\sigma\pi s)$-reflecting. Conversely, if an ordinal is $(\sigma\wedge \pi^kc^m_1s)$-reflecting, then, by the equation above, it is
\[\text{$\sigma\wedge \pi^k\sigma\pi(\sigma\pi^m\sigma\pi s\wedge \pi^m\sigma\pi s)$-reflecting},\]
and thus 
\[\text{$\sigma\wedge \pi^k\sigma\pi\sigma\pi^m\sigma\pi s$-reflecting}.\]
If so, then by Theorem \ref{secondmainwords}, it is
\[\text{$\sigma\wedge \pi^k\sigma\pi s$-reflecting},\]
as desired.
\endproof
\fi
\begin{lemma}\label{LemmaBelowOmegaCubed2}
For every $k \in\mathbb{N}$, every $n_0,\hdots, n_k \in\mathbb{N}$, and every reflection pattern $s$, every $(\sigma\wedge \pi^{k+1}s)$-reflecting ordinal is $(\pi^{k+1} \sigma\pi c^0_{n_0}c^1_{n_1}\dots c^k_{n_k}s)$-reflecting.
\end{lemma}
\proof
This follows from applying Corollary \ref{LemmaBelowOmegaCubedCorollary} repeatedly.
\endproof

\begin{definition}
A reflection pattern is in \emph{$2$-normal form} if it is of the form
\[\pi^{n}c^0_{n_0}c^1_{n_1}\dots c^k_{n_k},\]
for some natural numbers $k$, $n_0,n_1,\hdots, n_k$. If $w$ is the reflection pattern above, we define \index{$2$-normal form}
\[o(w) = \omega^{k+1}\cdot n_k + \omega^{k}\cdot n_{k-1} + \dots + \omega \cdot n_0 + n.\]
\end{definition}
For now, we shall simply refer to patterns in $2$-normal form as being in \emph{normal form}.
We shall see that patterns in normal form have very nice properties.

\begin{lemma}\label{LemmaSigmaSigmaUniversal}
Suppose $s$ is a reflection pattern in normal form. Then, every $\sigma\sigma$-reflecting ordinal is $\sigma(s)$-reflecting. \index{Contraction Lemma}
\end{lemma}
\proof
By Theorem \ref{firstmainwords}, any such ordinal is $\sigma(\sigma\wedge \pi^k)$-reflecting for any $k$. The lemma now follows from Lemma \ref{LemmaBelowOmegaCubed2}.
\endproof

We also have the following ``contraction'' lemma, which will be crucial:
\begin{lemma}\label{LemmaContractionWords}
Suppose $s$ is a reflection pattern. Then, every $\sigma\pi\sigma\pi s$-reflecting ordinal is $\sigma\pi s$-reflecting.
\end{lemma}
\proof
Suppose $\alpha$ is $\sigma\pi\sigma\pi s$-reflecting. By an argument as in Lemma \ref{Propswp<sseasy}, applying Lemma \ref{Propswp<ss}, for every $\Sigma^1_1$-sentence $\phi$ satisfied by $L_\alpha$, one can find some $\pi\sigma\pi s$-reflecting $\beta <\alpha$ such that $L_\beta\models\phi$ and $\beta$ is not $\Sigma^1_1$-reflecting.  By Theorem \ref{secondmainwords}, $\beta$ is $\pi \pi s$-reflecting. By Lemma \ref{Propswp<sseasy}, $\beta$ is $\pi s$-reflecting, as desired.
\endproof

\begin{lemma}\label{LemmaSSNFPi}
Suppose $\pi t$ and $\pi s$ are reflection patterns in normal form such that $o(\pi t) \leq o(\pi s)$. Then, every $(\pi s)$-reflecting ordinal is either $(\pi t)$-reflecting or $\sigma$-reflecting.
\end{lemma}
\proof
Let 
\[\pi t = \pi^{n}c^0_{n_0}c^1_{n_1}\dots c^k_{n_k},\]
and 
\[\pi s = \pi^{m}c^0_{m_0}c^1_{m_1}\dots c^l_{m_l},\]
where $n$ and $m$ are nonzero. Without loss of generality, we assume that $n_k$ and $m_l$ are also nonzero. It follows that $k \leq l$.
Suppose that $o(\pi t) < o(\pi s)$. It will be convenient, for illustrative purposes, to consider the case that $k < l$ first. If so, it suffices to show that every $\pi(\sigma\wedge\pi^{l})$-reflecting ordinal which is not $\sigma$-reflecting is $(\pi t)$-reflecting, for then the result follows from Theorem \ref{secondmainwords}. By Lemma \ref{LemmaBelowOmegaCubed2}, every $\pi(\sigma\wedge\pi^{l})$-reflecting ordinal is
\begin{align*}
\text{$\pi(\sigma\wedge \pi^l \sigma\pi c^0_{n_0}c^1_{n_1}\dots c^k_{n_k})$-reflecting.}
\end{align*}
By Lemma \ref{Propswp<ss}, every such ordinal is 
\[\text{$\pi(\sigma \pi^l c^0_{n_0}c^1_{n_1}\dots c^k_{n_k})$-reflecting,}\] 
and, by Theorem \ref{firstmainwords} or Lemma \ref{Propswp<sseasy}, according as $n \leq l$ or $l \leq n$, it is
\[\text{$\pi(\sigma \pi^n c^0_{n_0}c^1_{n_1}\dots c^k_{n_k})$-reflecting,}\] 
so that, if it is not $\sigma$-reflecting, then it is 
\[\text{$\pi(\pi^n c^0_{n_0}c^1_{n_1}\dots c^k_{n_k})$-reflecting,}\] 
by Theorem \ref{secondmainwords}, i.e., $\pi(t)$-reflecting.

The general case is similar: let $i\leq k$ be greatest such that $n_i < m_i$ and notice that
\[c_{m_i}^i = c_{m_i-n_i}^ic_{n_i}^i.\]
Thus, 
\[\pi t = \pi^{n}c^0_{n_0}c^1_{n_1}\dots c^{i-1}_{n_{i-1}}t',\]
where $t' = c^{i}_{n_{i}}\dots c^k_{n_k}$; and
\[\pi s = \pi^{m}c^0_{m_0}c^1_{m_1}\dots c_{m_i-n_i}^it'.\]
It suffices to show that every $\pi(\sigma\wedge\pi^i\sigma\pi t')$-reflecting ordinal which is not $\sigma$-reflecting is $(\pi t)$-reflecting, for then the result follows from Theorem \ref{secondmainwords}. 
Lemma \ref{LemmaBelowOmegaCubed2} (with $\sigma\pi t'$ being the $s$ in the statement) shows that every such ordinal is
\begin{align*}
\text{$\pi(\sigma\wedge \pi^i\sigma\pi c^0_{n_0}c^1_{n_1}\dots c^{i-1}_{n_{i-1}}\sigma\pi t')$-reflecting.}
\end{align*}
As before, by Lemma \ref{Propswp<ss}, every such ordinal is  
\begin{align*}
\text{$\pi(\sigma\pi^i\sigma\pi c^0_{n_0}c^1_{n_1}\dots c^{i-1}_{n_{i-1}}\sigma\pi t')$-reflecting.}
\end{align*}
By Theorem \ref{firstmainwords} and Lemma \ref{Propswp<sseasy}, every such ordinal is
\begin{align*}
\text{$\pi(\sigma\pi^n\sigma\pi c^0_{n_0}c^1_{n_1}\dots c^{i-1}_{n_{i-1}}\sigma\pi t')$-reflecting.}
\end{align*}
By contraction (Lemma \ref{LemmaContractionWords}), it is
\begin{align*}
\text{$\pi(\sigma\pi^n\sigma\pi c^0_{n_0}c^1_{n_1}\dots c^{i-1}_{n_{i-1}}t')$-reflecting,}
\end{align*}
and by an argument like the one for Lemma \ref{LemmaContractionWords}, it is 
\begin{align*}
\text{$\pi(\sigma\pi^n c^0_{n_0}c^1_{n_1}\dots c^{i-1}_{n_{i-1}}t')$-reflecting,}
\end{align*}
so that if it is not $\sigma$-reflecting, then it is
\begin{align*}
\text{$\pi(\pi^n c^0_{n_0}c^1_{n_1}\dots c^{i-1}_{n_{i-1}}t')$-reflecting,}
\end{align*}
by Theorem \ref{secondmainwords}, as desired.
\endproof

\begin{lemma}\label{LemmaSSNFSigma}
Suppose $t$ is a reflection pattern in normal form. Then $\sigma\wedge t$ is equivalent to a reflection pattern in normal form.
\end{lemma}
\proof
Put $t = \pi^{k}c^0_{m_0}c^1_{m_1}\dots c^l_{m_l}$. The lemma is immediate unless $k \neq 0$ and there is some least $i \leq l$ such that $m_i \neq 0$. Thus, 
\begin{align*}
t 
&= \pi^{k}c^i_{m_i}c^{i+1}_{m_{i+1}}\dots c^l_{m_l}\\
&= \pi^{k}(\sigma \wedge \pi^i \sigma\pi c^i_{m_i-1}c^{i+1}_{m_{i+1}}\dots c^l_{m_l}).
\end{align*}
Let
\[s = c^{k+i}_{m_{k+i}}c^{k+ i+1}_{m_{k+i+1}}\dots c^l_{m_l},\]
so that 
\[t = \pi^{k}(\sigma \wedge \pi^i \sigma\pi c^i_{m_i-1}c^{i+1}_{m_{i+1}}\dots c^{k+i-1}_{m_{k + i-1}}s).\] 
If all the indicated $m_j$ are zero, then the result follows easily; otherwise, by Lemma \ref{Propswp<ss} and Lemma \ref{LemmaContractionWords},
\[t = \pi^{k}(\sigma \wedge \pi^i \sigma\pi c^i_{m_i-1}c^{i+1}_{m_{i+1}}\dots c^{k + i-1}_{m_{k + i-1}}\sigma\pi s).\] 
By Lemma \ref{Propswp<ss},
\begin{align*}
\sigma\wedge t 
&= \sigma \wedge \pi^{k}(\sigma \wedge \pi^i \sigma\pi c^i_{m_i-1}c^{i+1}_{m_{i+1}}\dots c^{k+i-1}_{m_{k+i-1}}\sigma\pi s)\\
&= \sigma \wedge \pi^{k}(\pi^i \sigma\pi c^i_{m_i-1}c^{i+1}_{m_{i+1}}\dots c^{k+i-1}_{m_{k+i-1}}\sigma\pi s)\\
&= \sigma \wedge \pi^{k+i}\sigma\pi c^i_{m_i-1}c^{i+1}_{m_{i+1}}\dots c^{k+i-1}_{m_{k+i-1}}\sigma\pi s)
\end{align*}
By Lemma \ref{LemmaBelowOmegaCubed2} on the one hand and Lemma \ref{Propswp<sseasy} and Lemma \ref{LemmaContractionWords} on the other,
\[\sigma \wedge \pi^{k+i}\sigma\pi s \equiv \sigma \wedge t,\]
but the reflection pattern on the left-hand side is readily seen to be equivalent to one in normal form.
\endproof

\begin{theorem}\label{TheoremSSNormalForm}
Let $s$ be a reflection pattern in which the string $\sigma\sigma$ does not occur. Then, it is equivalent to a reflection pattern in $2$-normal form.
\end{theorem}
\proof
This is proved by induction on the construction of $s$. Clearly, if $s$ is equivalent to a reflection pattern in normal form, then so too is $\pi s$. We only need consider the string $\sigma s$ in the case that $s$ is of the form
\[\pi^{m}c^0_{m_0}c^1_{m_1}\dots c^l_{m_l},\]
where $m$ is nonzero.
Then, it is easy to see that
\begin{align*}
\sigma\pi^{m}c^0_{m_0}c^1_{m_1}\dots c^l_{m_l} 
&= \sigma\pi c^0_{m_0}c^1_{m_1}\dots c^l_{m_l}\\
&= c^0_{m_0+1}c^1_{m_1}\dots c^l_{m_l}.
\end{align*}
Now, let $s$ be as above, and let
\[t = \pi^{n}c^0_{n_0}c^1_{n_1}\dots c^l_{n_l}.\]
We need to show that $s\wedge t$ is equivalent to a reflection pattern in normal form.
If both $n$ and $m$ are nonzero, then the result follows from Lemma \ref{LemmaSSNFPi}; so suppose that one of $m$ and $n$ is zero, so that
\[s\wedge t = \sigma \wedge s \wedge t.\]
Write $n = n_{-1}$ and $m = m_{-1}$ and 
let $i$ and $j$ be least such that $n_i$ and $m_j$ are nonzero, respectively. There are four cases to consider. The first one is that in which both $i$ and $j$ are equal to $0$. Then, there are reflection patterns $u$ and $v$, both in normal form, such that
\begin{align*}
t = \sigma\pi u
\end{align*}
and 
\begin{align*}
s = \sigma \pi v.
\end{align*}
Suppose without loss of generality that $o(\pi u) < o(\pi v)$. 
By Lemma \ref{LemmaSSNFPi}, every $(\pi v)$-reflecting ordinal which is not $\sigma$-reflecting is $(\pi u)$-reflecting. Thus, every $s$-reflecting ordinal is either $t$-reflecting or $\sigma\sigma$-reflecting, in which case it is also $t$-reflecting by Lemma \ref{LemmaSigmaSigmaUniversal}. 

The second case is that in which $i = 0$ but $j \neq 0$. Then, there are reflection patterns $u$ and $v$, both in normal form, such that
\begin{align*}
t = \sigma\pi u
\end{align*}
and 
\begin{align*}
\sigma\wedge s = \sigma \wedge \pi^{m_j}\sigma\pi v.
\end{align*}
By direct computation,
\begin{align*}
s \wedge t 
& \equiv \sigma\pi u \wedge \pi^{m_j}\sigma\pi v\\
& \equiv \sigma(\pi u\wedge \pi^{m_j}\sigma\pi v) \wedge \pi^{m_j}\sigma\pi v\\
& \equiv \sigma(\pi u\wedge \pi^{m_j}\sigma\pi v) \wedge \pi^{m_j}(\sigma(\pi u\wedge \pi^{m_j}\sigma\pi v)\wedge\sigma\pi v)\\
& \equiv \sigma\wedge \pi^{m_j}(\sigma(\pi u\wedge \pi^{m_j}\sigma\pi v)\wedge\sigma\pi v).
\end{align*}
Now, it is easily seen that
\begin{align*}
\pi^{m_j}\sigma(\pi u\wedge \pi^{m_j}\sigma\pi v)
&\equiv \pi^{m_j}\sigma(\pi u\wedge \pi^{m_j}\sigma\pi v \wedge \pi\sigma\pi v)\\
&\equiv \pi^{m_j}\Big(\sigma(\pi u\wedge \pi^{m_j}\sigma\pi v)\wedge \sigma\pi\sigma\pi v \Big) \\
&\equiv \pi^{m_j}\Big(\sigma(\pi u\wedge \pi^{m_j}\sigma\pi v)\wedge \sigma\pi v \Big),
\end{align*}
where the last equivalence follows from Lemma \ref{LemmaContractionWords},
and so
\[s\wedge t \equiv \sigma \wedge \pi^{m_j}\sigma(\pi u\wedge \pi^{m_j}\sigma\pi v).\]
Since each of $\pi u$ and $\pi^{m_j}\sigma\pi v$ is a reflection pattern in normal form, Lemma \ref{LemmaSSNFPi} implies that their conjunction is equivalent to one of $\pi u$ and $\pi^{m_j}\sigma\pi v$ in this context.  Let us denote this conjunct by $w$. Then, by an argument as before,
\begin{align*}
s \wedge t
&\equiv \sigma \wedge \pi^{m_j}\sigma w.
\end{align*}
Since $w$ is in normal form and of the form $\pi w'$, so too is $\pi^{m_j}\sigma w$, so the result follows from Lemma \ref{LemmaSSNFSigma}. The case in which $i \neq 0$ and $j = 0$ is analogous.

The remaining case is that in which both $i$ and $j$ are nonzero. Suppose without loss of generality that $n_i < m_j$
By replacing $n_i$ and $m_j$ by larger numbers if necessary (this might need to be done in the case $i = -1$ and $m_0 = 0$---a situation similar to the one in Lemma \ref{LemmaSSNFSigma}) we may assume that there are reflection patterns $u$ and $v$ such that
\begin{align*}
\sigma\wedge t = \sigma\wedge \pi^{n_i}\sigma \pi u
\end{align*}
and 
\begin{align*}
\sigma\wedge s = \sigma\wedge \pi^{m_j}\sigma \pi v.
\end{align*}
Then, we have
\begin{align*}
\sigma\wedge s \wedge t 
&\equiv\sigma\wedge \pi^{n_i}\sigma\pi u\wedge\pi^{m_j}\sigma\pi v\\
&\equiv\sigma\big(\pi^{n_i}\sigma\pi u\wedge\pi^{m_j}\sigma\pi v\big)\wedge \pi^{n_i}\sigma\pi u\wedge\pi^{m_j}\sigma\pi v\\
&\equiv\sigma\wedge \pi^{n_i}\sigma\pi u\wedge\pi^{m_j}\Big(\sigma\big(\pi^{n_i}\sigma\pi u\wedge\pi^{m_j}\sigma\pi v\big)\wedge \sigma\pi v\Big).
\end{align*}
By contraction,
\begin{align*}
\sigma \wedge \pi^{m_j}\sigma\pi^{n_i}\sigma\pi u
&\equiv \sigma \wedge \pi^{m_j}\sigma\pi\sigma\pi u\\
&\equiv \sigma \wedge \pi^{m_j}\sigma\pi u,
\end{align*}
so, because $n_i<m_j$, it follows that every ordinal which is 
\[\text{$(\sigma \wedge \pi^{m_j}\sigma\pi^{n_i}\sigma\pi u)$-reflecting}\] 
is also  
\[\text{$(\sigma \wedge \pi^{n_i}\sigma\pi u)$-reflecting.}\]
Similarly, we have
\begin{align*}
\sigma \wedge \pi^{m_j}\sigma\pi^{m_j}\sigma\pi v
&\equiv \sigma \wedge \pi^{m_j}\Big(\sigma\pi^{m_j}\sigma\pi v\wedge \sigma\pi\sigma\pi v\Big)\\
&\equiv \sigma \wedge \pi^{m_j}\Big(\sigma\pi^{m_j}\sigma\pi v\wedge \sigma\pi v\Big)
\end{align*}
and so every ordinal which is $(\sigma \wedge \pi^{m_j}\sigma\pi^{m_j}\sigma\pi v)$-reflecting is also 
\[\text{$\Big(\sigma \wedge \pi^{m_j}\big(\sigma\pi^{m_j}\sigma\pi v\wedge \sigma\pi v\big)\Big)$-reflecting.}\]
From these two observations, we see that
\begin{align*}
\sigma\wedge s \wedge t
&\equiv
\sigma\wedge \pi^{n_i}\sigma\pi u\wedge\pi^{m_j}\Big(\sigma\big(\pi^{n_i}\sigma\pi u\wedge\pi^{m_j}\sigma\pi v\big)\wedge \sigma\pi v\Big)\\
&\equiv \sigma\wedge \pi^{m_j}\sigma\big(\pi^{n_i}\sigma\pi u\wedge\pi^{m_j}\sigma\pi v\big).
\end{align*}
Since both $m_j$ and $n_i$ are nonzero, Lemma \ref{LemmaSSNFPi} implies that in this context the conjunction 
\[\pi^{n_i}\sigma\pi u\wedge\pi^{m_j}\sigma\pi v\]
is equivalent to one of $\pi^{n_i}\sigma\pi u$ or $\pi^{m_j}\sigma\pi v$. Denote this conjunct by $w$. Then, $w$ is in normal form. Since $w$ begins with the symbol $\pi$, $\pi^{m_j}\sigma w$ is also equivalent to a pattern in normal form. The result then follows from Lemma \ref{LemmaSSNFSigma}. This proves the theorem.
\endproof

\begin{theorem}
$\sigma\sigma$ has order-type $\omega^\omega$ in the reflection order.
\end{theorem}
\proof
By Theorem \ref{TheoremSSNormalForm}, every reflection pattern in which the string $\sigma\sigma$ does not occur is equivalent to one in $\sigma\sigma$-normal form. By Lemma \ref{LemmaSigmaSigmaUniversal}, this implies that $\sigma\sigma$ is strictly bigger than each reflection pattern in which the string $\sigma\sigma$ does not occur. Conversely, if a reflection pattern does contain the string $\sigma\sigma$, then naturally, it cannot be strictly smaller than $\sigma\sigma$. An easy induction using Lemma \ref{LemmaBelowOmegaCubed2} shows that, for reflection patterns $u$ and $v$ in $2$-normal form, $u < v$ if, and only if, $o(u) < o(v)$, so the result follows.
\endproof
\section{Linear Patterns}\label{SectConjunctionFree}
Our first result in this section concerns the length of the linear reflection order; its proof induces a simple algorithm for comparing two arbitrary linear reflection patterns.

\begin{theorem}
The length of the linear fragment of the reflection order is $\omega^\omega$.
\end{theorem}
\proof
Let us employ the simplified notation from the previous section.
Recursively, we assign ordinals to reflection patterns without conjunction: we assign the ordinal $\omega^n$ to the pattern
\[\sigma^n.\]
In particular, the ordinal $1$ is assigned to the empty pattern. If $s$ and $t$ are patterns to which ordinals $o(s)$ and $o(t)$ have been assigned, we assign the ordinal 
\[o(s) + o(t)\]
to the pattern
\[t \pi s.\]
Note that if $n < m$, then, on the one hand,
\[\alpha + \omega^n + \omega^m = \alpha + \omega^m,\]
while, on the other,
\begin{align*}
\sigma^{m} \pi s
&\equiv \sigma^{m-n}\sigma^n \pi s\\
&\equiv \sigma^{m-n}(\sigma^n \pi s \wedge \pi) & \text{by Lemma \ref{Propswp<ss}\eqref{Propswp<ss3}}\\
&\equiv \sigma^{m-n}(\sigma^n\pi s \wedge \pi\sigma^n \pi s) & \text{by Lemma \ref{Propswp<ss}\eqref{Propswp<ss1}}\\
&\equiv \sigma^{m-n}\Big(\sigma^n(\pi s \wedge \pi\sigma^n \pi s) \wedge \pi\sigma^n\pi s\Big) &  \text{by Lemma \ref{Propswp<ss}\eqref{Propswp<ss2}}\\
&\equiv \sigma^m(\pi s \wedge \pi\sigma^n \pi s) \\
&\equiv \sigma^m \pi s \wedge \sigma^m\pi\sigma^n \pi s.
\end{align*}
Thus, every $\sigma^m \pi s$-reflecting ordinal is also $\sigma^m \pi \sigma^n \pi s$-reflecting when $n < m$; the converse is also true, by Lemma \ref{Propswp<sseasy} and Theorem \ref{secondmainwords} (cf. the argument of Lemma \ref{LemmaContractionWords} on p. \pageref{LemmaContractionWords}). It follows that this assignment of ordinals is well-defined. By Lemma \ref{Propswp<sseasy} and Theorem \ref{firstmainwords}, we have
\[\sigma\pi^n \sigma \equiv \sigma\pi^m\sigma\]
for any pair of nonzero numbers $n$ and $m$, so every conjunction-free reflection pattern is equivalent to one to which an ordinal has been assigned. It should be clear by now that, for conjunction-free reflection patterns, $s < t$ if, and only if, $o(s) < o(t)$, which completes the proof.
\endproof

The second result of this section is that the linear patterns are cofinal in the reflection order.

\begin{theorem}\label{theoremcofinalsequence}
The sequence $\{(\sigma^1_1)^n:n\in\mathbb{N}\}$ is cofinal in the reflection order.
\end{theorem}
\proof
To prove the theorem, we shall prove by induction on the construction of a reflection pattern $s$ that if $s$ contains no occurrence of $\sigma^{n+1}$, then every $\sigma(\sigma^n \wedge t)$-reflecting ordinal is also $\sigma(\sigma^n\wedge t \wedge s)$-reflecting, for every reflection pattern $t$ (cf. Lemma \ref{LemmaSigmaSigmaUniversal} on p. \pageref{LemmaSigmaSigmaUniversal}). We may assume that $n \neq 0$, for otherwise the conclusion follows from Theorem \ref{firstmainwords}.
Note that the case that $s$ is a conjunction is immediate from the induction hypothesis and the case that $s$ is of the form $\pi s'$ is also immediate from the induction hypothesis and Theorem \ref{firstmainwords}, thus, we suppose that $s$ is of the form $\sigma s'$. The pattern $s'$ might be a conjunction, say,
\[s' = s_0\wedge s_1 \wedge \dots\wedge s_k \wedge r_0 \wedge r_1 \wedge \dots \wedge r_l,\]
where each $s_i$ is of the form $\sigma^{m_i}s_i'$ for some $m_i < n$ and some $s_i'$ which is a conjunction of patterns of the form $\pi s^*$, and each $r_i$ is of the form $\pi r_i'$. Instead of proving that every $\sigma(\sigma^n\wedge t)$-reflecting ordinal is $\sigma(\sigma^n\wedge t \wedge \sigma s')$-reflecting, we shall prove the stronger fact that it is
\[\sigma\Big(\sigma^n\wedge t \wedge \sigma\big(\sigma^{n-1}(s_0'\wedge s_1'\wedge\dots \wedge s_k')\wedge r_0\wedge r_1 \wedge \dots \wedge r_l\big)\Big)\text{-reflecting.}\]
By the induction hypothesis applied to $s'$, and the fact that each of $s'_i$ and
\[r :=  r_0 \wedge r_1 \wedge \dots \wedge r_l\]
is a conjunction of patterns of the form $\pi s^*$, we obtain:
\begin{align*}
\sigma\Big(\sigma^n\wedge t\Big)
&\equiv \sigma\Big(\sigma^n\wedge t\wedge s_0' \wedge \dots\wedge s_k' \wedge r\Big)\\
&\equiv \sigma\Big(\sigma^n\wedge t\wedge \sigma^n\wedge s_0' \wedge \dots\wedge s_k' \wedge r\Big)\\
&\equiv \sigma\Big(\sigma^n\wedge t\wedge \sigma\big(\sigma^{n-1}\wedge s_0' \wedge \dots\wedge s_k' \wedge r\big)\Big)\\
&\equiv \sigma\Big(\sigma^n\wedge t\wedge \sigma\big(\sigma^{n-1}(s_0' \wedge \dots\wedge s_k') \wedge r\big)\Big),
\end{align*}
where the last two equivalences follow from Lemma \ref{Propswp<ss}. This completes the proof of the theorem.
\endproof

\section{Concluding Remarks}\label{SectConcludingPatterns}
We have not given any bounds on the length of the reflection order. Let us say something about this:
\begin{proposition} \label{propBoundRO}
Let $s$ be a reflection pattern and suppose $\alpha$ is a countable, locally countable ordinal such that
\[L_\alpha\prec_1 L_{\alpha^++1}.\]
Then, $\alpha$ is $s$-reflecting.
\end{proposition}
\proof
Since 
\[L_\alpha\prec_1 L_{\alpha^+},\]
it follows that $\alpha$ is $\pi^1_1$-reflecting. By Gostanian's theorem \cite{Go79} mentioned in the introduction, $\alpha$ is $\sigma^1_1$-reflecting. Inductively, suppose it is
\[(\sigma^1_1)^n\text{-reflecting}\]
and let $\psi$ be a $\Sigma^1_1$ sentence such that
\[L_\alpha\models \psi.\]
Choose a $\Pi_1$ sentence $\psi^*$ such that for all admissible $\beta$ containing all relevent parameters,
\[L_\beta\models\psi\]
if, and only if,
\[L_{\beta^+}\models\psi^*(L_\beta),\]
so that, in particular,
\[L_{\alpha^+}\models\psi^*(L_\alpha).\]
Then, from the point of view of $L_{\alpha^++1}$, there are admissible sets $L_\alpha$ and $L_{\alpha^+}$ such that
\begin{enumerate}
\item $\alpha$ is $(\sigma^1_1)^n$-reflecting, i.e., for every $\Pi_1$ sentence $\phi^*$, if $L_{\alpha^+}\models\phi^*(L_\alpha)$, then there is a $(\sigma^1_1)^{n-1}$-reflecting $\beta<\alpha$ such that $L_{\beta^+}\models\phi^*(L_\beta)$. (The quantification over $\phi^*$ is bounded.)
\item $L_{\alpha^+}\models\psi^*(L_\alpha)$.
\end{enumerate}
Thus, by stability, there are admissible sets $L_\beta$ and $L_{\beta^+}$ in $L_\alpha$ such that $\beta$ is $(\sigma^1_1)^n$-reflecting and $L_{\beta^+}\models\psi^*(L_\beta)$. Hence, $\alpha$ is $(\sigma^1_1)^{n+1}$-reflecting. Therefore, a simple induction on the construction of a reflection pattern $s$, using the proofs of Theorem \ref{secondmainwords} and Theorem \ref{theoremcofinalsequence} shows that $\alpha$ is $s$-reflecting.
\endproof

The length of the reflection order is thus at most the least ordinal $\alpha$ such that
\[L_\alpha\prec_1 L_{\alpha+1}.\]
Moreover, surely each inequality between reflection patterns is provable in any theory that proves the existence of the corresponding ordinals. This suggests strongly that the length of the reflection order is smaller than the proof-theoretic ordinal of the subsystem $\PI^1_2$-CA$_0$ of analysis and in fact smaller than the ordinal described in Rathjen \cite{Ra05a}, though we do not have a proof of this. The reader may consult \cite{Ag} for an example of a chain of length $\varepsilon_0$ in the reflection order.

An interesting question is that of the structure of the ``higher'' reflection order, defined in terms of iterated $\Pi^1_n$ and $\Sigma^1_n$-reflection and conjunctions. The situation there is very different and involves set-theoretic considerations; it will be the subject of a forthcoming article.

\bibliographystyle{abbrv}
\bibliography{References}

\begin{thebibliography}{10}

\bibitem{Aa74}
S.~Aanderaa.
\newblock {Inductive Definitions and Their Closure Ordinals}.
\newblock In J.~E. Fenstad and P.~G. Hinman, editors, {\em Generalized
  Recursion Theory}, pages 207--220. 1974.

\bibitem{AbSa76}
F.~G. Abramson and G.~E. Sacks.
\newblock {Uncountable Gandy Ordinals}.
\newblock {\em J. London Math. Soc.}, 14(2):387--392, 1976.

\bibitem{AcRi74}
P.~Aczel and W.~Richter.
\newblock {Inductive Definitions and Reflecting Properties of Admissible
  Ordinals}.
\newblock In J.~E. Fenstad and P.~G. Hinman, editors, {\em Generalized
  Recursion Theory}, pages 301--381. 1974.

\bibitem{Ag}
J.~P. Aguilera.
\newblock {\em {Between the Finite and the Infinite}}.
\newblock 2019.
\newblock Ph.D. Thesis. Vienna University of Technology.

\bibitem{Ba75}
J.~Barwise.
\newblock {\em Admissible Sets and Structures}.
\newblock Perspectives in Mathematical Logic. Springer-Verlag, Berlin, 1975.

\bibitem{BGM}
K.~J. Barwise, R.~Gandy, and Y.~N. Moschovakis.
\newblock {The Next Admissible Set}.
\newblock {\em J. Symbolic Logic}, 36:108--120, 1971.

\bibitem{Ce74}
D.~Cenzer.
\newblock {Ordinal Recursion and Inductive Definitions}.
\newblock In J.~E. Fenstad and P.~G. Hinman, editors, {\em Generalized
  Recursion Theory}, pages 221--264. 1974.

\bibitem{Go79}
R.~Gostanian.
\newblock {The Next Admissible Ordinal}.
\newblock {\em Ann. Math. Logic}, 17:171--203, 1979.

\bibitem{GoHr79}
R.~Gostanian and K.~Hrbacek.
\newblock {A New Proof That $\pi^1_1<\sigma^1_1$}.
\newblock {\em Zeitsch. f. math. Logik und Grundlagen d. Math.}, 25:407--408,
  1979.

\bibitem{Ra05a}
M.~Rathjen.
\newblock {An Ordinal Analysis of Parameter-Free $\Pi^1_2$-Comprehension}.
\newblock {\em Arch. Math. Logic}, 44:263--362, 2005.

\bibitem{Si78}
S.~G. Simpson.
\newblock {Short Course on Admissible Recursion Theory}.
\newblock In J.~E. Fenstad, R.~O. Gandy, and G.~E. Sacks, editors, {\em
  Generalized Recursion Theory, II}, pages 355--390. 1978.

\end{thebibliography}

\end{document}